\documentclass[11pt]{article}
\usepackage{graphicx}
\usepackage{amsmath,amssymb,amsthm,amsfonts}
\usepackage{amsmath}
\allowdisplaybreaks[4]
\usepackage{color}
\usepackage{amssymb}
\usepackage{cite}
\usepackage{appendix}
\newtheorem{thm}{Theorem}[section]

\newtheorem{lem}[thm]{Lemma}

\theoremstyle{definition}

\theoremstyle{remark}

\numberwithin{equation}{section}

\DeclareMathSymbol{\C}{\mathalpha}{AMSb}{"43}

\newcommand{\eps}{\varepsilon}

\newcommand{\R}{{\mathbb{R}}}

\newcommand{\inte}{\int_{\mathbb{R}^N}}

\newcommand{\bsub}{\begin{subequations}}
\newcommand{\esub}{\end{subequations}$\!$}

\begin{document}
\title{Limiting Behavior of Constraint Minimizers for Inhomogeneous Fractional Schr\"{o}dinger Equations}

\author{Hongfei Zhang\thanks{School of Mathematics and Statistics, Central China Normal University, P.O. Box 71010, Wuhan 430079,
P. R. China.  Email: \texttt{hfzhang@mails.ccnu.edu.cn}.}
\, and\, Shu Zhang\thanks{School of Mathematics and Statistics, Central China Normal University, P.O. Box 71010, Wuhan 430079,
P. R. China.  Email: \texttt{szhangmath@mails.ccnu.edu.cn}.}
}

\date{}

\smallbreak \maketitle

\begin{abstract}
We study $L^2$-normalized solutions of the following inhomogeneous fractional Schr\"{o}dinger equation
\begin{equation*}
(-\Delta)^{s} u(x)+V(x)u(x)-a|x|^{-b}|u|^{2\beta^2}u(x)=\mu u(x)\ \ \mbox{in}\ \ \R^{N}.
\end{equation*} 
Here $s\in(\frac{1}{2},1)$, $N>2s$, $a>0$, $0<b<\min\{\frac{N}{2},1\}$, $\beta=\sqrt{\frac{2s-b}{N}}$ and $V(x)\geq 0$ is an external
potential. We get $L^2$-normalized solutions of the above equation by solving the associated
constrained minimization problem. We prove that there exists a critical value $a^*>0$ such that minimizers exist for $0<a<a^*$, and minimizers do
not exist for any $a>a^*$. In the case of $a=a^*$, one can obtain the classification results of the existence and non-existence for constraint minimizers, which are depended strongly on the value of $V(0)$. For $V(0)=0$, the limiting behavior of nonnegative minimizers is also analyzed when $a$ tend to $a^*$ from below. 
\end{abstract}

\vskip 0.05truein


\section{Introduction}
We consider the following inhomogeneous fractional Schr\"{o}dinger equation (NLS)
\begin{equation}\label{F}
(-\Delta)^{s} u(x)+V(x)u(x)-a|x|^{-b}|u|^{2\beta^2}u(x)=\mu u(x)\ \ \mbox{in}\ \ \R^{N},
\end{equation}
where $s\in(\frac{1}{2},1)$, $N>2s$, $a>0$, $0<b<\min\{\frac{N}{2},1\}$ and $\beta=\sqrt{\frac{2s-b}{N}}$ denotes the $L^2$-critical exponent, $V(x)$ and  $\mu\in\R$ describe the trapping potential and the chemical potential, respectively. Note from \cite{D2} that $(-\Delta)^{s} u$ is a fractional Laplace operator satisfying \begin{equation}\label{01}
(-\Delta)^su(x)=-\frac{C(N,s)}{2}\int_{\R^N}\frac{u(x+y)+u(x-y)-2u(x)}{|y|^{N+2s}}dy\ \ \mbox{for any}\ \ x\in \R^N.
\end{equation}
Here $C(N,s)$ is a positive constant depending on $N$ and $s$, given by
\begin{equation}\label{02}
C(N,s)=\Big(\int_{\R^N}\frac{1-cos \xi_1}{|\xi|^{N+2s}}d\xi\Big)^{-1}.
\end{equation}

Equation \eqref{F} arises in various physical contexts like nonlinear optics, plasma physics and Bose-Einstein condensates (BECs) (cf.\cite{AG,Liu,Za}). The nonlocal operator $(-\Delta)^s$ can be seen as the infinitesimal generators of L$\acute{e}$vy stable diffusion processes \cite{A}, which appears in many physical models such as water models, quantum mechanics, the fractional Brownian motion \cite{D, RS} and so on.  The function $a|u|^{2\beta^2}|x|^{-b}$ in \eqref{F} stands for corrections to the nonlinear power-law response or some inhomogeneities in the medium. From the physical point of view, we assume that the trapping potential $V(x)$ satisfies	
\begin{itemize}
	\item[\rm($V_{1})$.] $0\leq V(x)\in {L^{\infty}_{loc}} (\R^{N})$  and $ \lim_{|x|\rightarrow\infty}V(x)=\infty$.
\end{itemize}

In this paper, we are concerned with the $L^2$-constraint variational problem for inhomogeneous fractional Schr\"{o}dinger equation \eqref{F}, which can be described by the following constrained minimization problem
\begin{equation}\label{e(a)}
e(a)=\inf _{\{u\in\mathcal{H}, \|u\|^{2}_{2}=1\}}E_{a}(u),\,\ a>0,
\end{equation}
where energy functional $E_{a}(u)$ is given by
\begin{equation}\label{Eau}
E_{a}(u):=\int_{\R^{N}}(|(-\Delta)^{\frac{s}{2}} u|^{2}+V(x)|u|^{2})dx-\frac{a}{1+\beta^2}\int_{\R^{N}}\frac{|u(x)|^{2+2\beta^2}}{|x|^{b}}dx,
\end{equation}
and
\begin{equation}\label{04}
\int_{\R^{N}}|(-\Delta)^{\frac{s}{2}} u|^{2}dx=\frac{C(N,s)} {2}\iint_{\R^{2N}}\frac{|u(x)-u(y)|^2}{|x-y|^{N+2s}}dxdy.
\end{equation}
The $\mathcal{H}$ in \eqref{e(a)} is defined as
\begin{equation*}
\mathcal{H}:=\Big\{u\in H^{s}(\R^{N}):\int_{\R^{N}}V(x)|u(x)|^{2}dx<\infty\Big\}
\end{equation*}
with associated norm
\begin{equation*}
\|u\|_{\mathcal{H}}=\Big\{\displaystyle\int_{\R^N}\Big(|(-\Delta)^{\frac{s}{2}} u(x)|^{2}+[1+V(x)]|u(x)|^{2}\Big)dx\Big\}^{\frac{1}{2}}.
\end{equation*}

When $s=1$, \eqref{e(a)} is the so-called $L^2$-constraint variational problem of Schr\"{o}dinger equations. In the last two decades, the existence and nonexistence of constrained minimizers and various quantitative properties of these equations have been investigated extensively, see \cite{AD,FG,SC,GR} and the references therein. More specifically, for $b=0$, (\ref{e(a)}) is a homogeneous problem. Guo, Zhang and their collaborators in \cite{GR,ZJ} proved that there is a critical constant $a_0>0$ such that minimizers of $e(a)$ exist if and only if $0<a<a_0.$  For different kinds of trapping potentials, the authors in \cite{GR,GZ} also analyzed the detailed limiting  behavior of non-negative minimizers as $a\nearrow a_0$. By constructing Pohozaev identities,  Guo, Lin, and Wei in \cite{GLW} obtained the uniqueness of positive minimizers for $e(a)$ as $a\nearrow a_0$. In the case of $b\neq0$,  (\ref{e(a)}) is an inhomogeneous problem which contains the nonlinear term $K(x)|u(x)|^{2+2\beta^2}$ with $K(x)=|x|^{-b}$. Due to the singularity of $|x|^{-b}$ at $x=0$, the mathematical analysis of the equation \eqref{F} however becomes more challenging, which has attracted a lot of attentions over the past few decades, see  \cite{AD,FG,SC} and the references therein. Especially, Genoud in \cite{FG} studied the existence of minimizers for \eqref{e(a)} in the special case $V(x)\equiv 1$.

Recently, the $L^2$-constraint variational problems involving in a fractional Laplacian have been widely studied, see \cite{D1,DT,HL} and the references therein. For  $0<s<1$, $b=0$, $\beta^2=\frac{4s}{N}$, and $V(x)$ satisfies $(V_1)$ or
\begin{itemize}
\item[\rm($V_{2})$.] $V(x)\in C(\R^N),\ \ \ 0=\inf_{x\in\R^N }V(x)<\sup_{x\in\R^N }V(x)=\lim_{|x|\rightarrow\infty}V(x)<\infty$.
\end{itemize}
Applying the classical concentration compactness principle \cite{L1, L2}, He and Long in \cite{HL} proved the existence and nonexistence of minimizers for $e(a)$. 
Subsequently, for $0<s<1$, $b=0$ and $0<\beta^2\leq\frac{4s}{N}$, making full use of fractional Gagliardo-Nirenberg-Sobolev inequality \cite[Proposition 3.1]{FLS}, the authors in \cite{D1} and \cite{DT} gave a complete classification of the existence and nonexistence of minimizers for $e(a)$ with periodic potentials and  trapping potentials, respectively. Furthermore, in the case of mass critical, under a suitable assumption of potential $V(x)$, they also analyzed the limiting behavior of minimizers when the mass tends to a critical value.

Motivated by above works, the main purpose of this paper is to investigate  the existence and nonexistence of minimizers for $e(a)$, and the limiting behavior of nonnegative minimizers. We first recall the following fractional Gagliardo-Nirenberg inequality.
\vskip 0.05truein

\noindent {\bf Theorem A (\cite[Theorem 2.2]{PZ}).} {\em
Let $s\in(\frac{1}{2},1)$, $N>2s$, $0<b<\min\{\frac{N}{2},1\}$ and $\beta^2=\frac{2s-b}{N}$. Then we have
\begin{equation}\label{GN}
\frac{a^*}{1+\beta^2}\int_{\R^{N}}|x|^{-b}|u(x)|^{2+2\beta^2}dx\leq \|(-\Delta)^{\frac{s}{2}} u(x)\|^{2}_{2}\|u(x)\|^{2\beta^2}_{2} \ \mbox{in}\ \ H^{s}(\R^{N}),
\end{equation}
where the optimal constant $a^*>0$ is given by
\begin{equation}\label{a8}
a^*:=\|Q\|^{2\beta^2}_{2}.
\end{equation}
Here $Q> 0$ is an optimizer of \eqref{GN} and solves the following equation
\begin{equation}\label{Eq}
(-\Delta)^{s} u+u-|x|^{-b}|u|^{2\beta^2}u=0\  \ \mbox{in}\ \ \R^{N}.
\end{equation}
Moreover, $Q>0$ satisfies
\begin{equation}\label{DS}
\int_{\R^{N}}|(-\Delta)^{\frac{s}{2}} Q|^{2}dx=\frac{1}{\beta^2}\int_{\R^{N}}|Q|^{2}dx=\frac{1}{1+\beta^2}\int_{\R^{N}}|x|^{-b}|Q|^{2+2\beta^2}dx.
\end{equation}
}

By making full use of the Gagliardo-Nirenberg inequality (\ref{GN}), the existence and nonexistence of minimizers for (\ref{e(a)}) can be stated as follows.

\begin{thm}\label{th1}
Let $s\in(\frac{1}{2},1)$, $N>2s$, $0<b<\min\{\frac{N}{2},1\}$, $\beta^2=\frac{2s-b}{N}$, $a^*=\|Q\|^{2\beta^2}_{2}$, and $V(x)$ satisfies $(V_{1})$. Then we have
\begin{enumerate}
\item If $0<a< a^{\ast}$, there exists at least one minimizer for (\ref{e(a)}).
\item If $a>a^{\ast}$, there is no minimizer for (\ref{e(a)}).
\end{enumerate}
Moreover, we also have $e(a)>0$ for $a<a^{\ast}$ and $e(a)=-\infty$ for $a>a^{\ast}$.
\end{thm}
There are several remarks on the proof of Theorem \ref{th1}. We first note that $(-\Delta)^s$ with $\frac{1}{2} < s < 1$ is a nonlocal operator, which causes extra difficulties. In particular, we need to estimate the Gagliardo (semi) norm of trial function, see \eqref{Gj} in Section 2. Inspired by \cite[Proposition 21]{SV}, we establish Lemma \ref{lem2.2} to circumvent this obstacle. Moreover,  different from the case $s = 1$ in \cite{LZ}, the function $Q(x)$ given by Theorem A is polynomially decay at infinity, see Lemma \ref{lem2.3} below.  We thus need to give more detailed analyses to establish the desired estimates of the trial function. Finally, due to the singularity of $|x|^{-b}$, we shall be carefully to verify \eqref{2.1}.

In view of Theorem \ref{th1}, a natural question is whether there exist minimizers for $e(a)$ when $a=a^*$.
In contrast to the homogeneous $b=0$,
our following results show that $e(a)$ may admit minimizers at the threshold $a^*$ for some trapping potentials $V(x)$ satisfying $(V_1)$.
Indeed, applying the Gagliardo-Nirenberg inequality, we deduce from the nonnegativity of $V(x)$ that $e(a^*)\geq 0$. On the other hand, taking the same test function as in \eqref{st} and letting $\tau\to\infty$, one can obtain that $e(a^*)\leq E_a(u_\tau)\leq V(0)$. Therefore, one can conclude from above that
$$0\leq e(a^*)\leq V(0)$$
and the result for the case of $a=a^*$ can be summarized as follows.
\begin{thm}\label{th2}
Let $s\in(\frac{1}{2},1)$, $N>2s$, $0<b<\min\{\frac{N}{2},1\}$, $\beta^2=\frac{2s-b}{N}$, $a^*=\|Q\|^{2\beta^2}_{2}$, and $V(x)$ satisfies $(V_{1})$. Then we have
\begin{enumerate}
\item If $0\leq e(a^*)<V(0)$, there exists at least one minimizer for $e(a^{\ast})$.
\item If  $V(0)=0$, there is no minimizer for $e(a^{\ast})$. Moreover, $\lim_{a\nearrow a^{\ast}}e(a)=e(a^{\ast})=0$.
\end{enumerate}
\end{thm}
Theorem \ref{th2} shows a fact that the existence and nonexistence of minimizers for $e(a^*)$ depend strongly on the value of $e(a^*)$ and $V(0)$. In brief, $V(0)$ can be regarded as a threshold of the energy $e(a^*)$ such that \eqref{e(a)} admits minimizers if $e(a^*)$ is less than $V(0)$, while \eqref{e(a)} does not admit any minimizer if $V(0)= e(a^*)=0$. In view of Theorem \ref{th2} (2), a simple analysis yields that minimizers must blow up in the sense of $\inte |(-\Delta)^{\frac{s}{2}} u_a|^2dx\to\infty$ as $a\nearrow a^*$. Therefore, we next focus on the limiting  behavior of minimizers $u_a$ as $a\nearrow a^*$. Note that $E_{a}(u_a)\geq E_{a}(|u_a|)$, then $|u_a|$ is also a minimizer of $e(a)$. Therefore, without loss of generality, one can restrict minimizers of $e(a)$ to nonnegative functions. We now establish the following theorem on limiting  behavior of nonnegative minimizers for $e(a)$ as $a\nearrow a^*$.

\begin{thm}\label{th3}
Let $s\in(\frac{1}{2},1)$, $N>2s$, $0<b<\min\{\frac{N}{2},1\}$,  $\beta^2=\frac{2s-b}{N}$ and $a^*=\|Q\|^{2\beta^2}_{2}$. Suppose $V(x)$ satisfies $(V_{1})$ with $V(0)=0$. For any $0<a<a^*$, assume $u_{a}$ is a nonnegative minimizer  of $e(a)$, define
\begin{equation}\label{jyp}
\varepsilon_{a}:=\Big(\int_{\R^{N}}|(-\Delta)^{\frac{s}{2}} u_{a}|^{2}dx\Big)^{-\frac{1}{2s}},
\end{equation}
and
\begin{equation}\label{jxpa}
w_{a}(x):=\varepsilon^{\frac{N}{2}}_{a}u_{a}(\varepsilon_{a}x).
\end{equation}
Then there exists a subsequence, still denoted by $\{w_a\}$, of $\{w_a\}$ such that
\begin{equation}\label{bp}
\lim\limits_{a\nearrow a^*}w_a=w_0\ \ \hbox{strongly in}\ \ H^s(\R^N),
\end{equation}
where $w_0$ is an optimizer of  the Gagliardo-Nirenberg inequality \eqref{GN}.
In addition, we have 
\begin{equation}\label{jbp}
\|w_a(x)\|_{L^{\infty}(\R^N)}\leq C\ \ \mbox{and}\ \ w_a(x)\to 0\ \ \mbox{as}\ \ |x|\to\infty\ \ \mbox{uniformly for}\ \ a\nearrow a^*,
\end{equation}
where $C>0$ is independent of $0<a<a^*$.
\end{thm}

There are several remarks on Theorem \ref{th3} in sequence. Firstly, since the  limit equation \eqref{Eq} has no translation invariance, which yields that $u_a$ must concentrate at the origin rather than any global minimum point of $V(x)$ as $a\nearrow a^*$. Secondly, it is unfortunately open whether the optimizer of \eqref{GN} is unique, so the convergence of theorem \ref{th3} can only holds in the sense of a subsequence. Moreover, different from \cite [Theorem 1.5] {LZ}, the De Giorgi-Nash-Moser type theorem of fractional Laplacian equations cannot be used directly to prove \eqref{jbp}. In order to overcome this difficulty, we first obtain $w_a\in L^{\infty}(\R^N)$ by Moser iteration \cite{M1}. Stimulated by \cite{AM,AV}, we then derive that $w_a$ tends to zero at infinity as $a\nearrow a^*$. Comparing with \cite{AM,AV}, there are some extra difficulties due to the singularity of $|x|^{-b}$, which is one of the innovations and difficulties of this paper.

This paper is organized as follows. In Section 2, we shall prove Theorems \ref{th1} and \ref{th2} on the existence and nonexistence of minimizers for $e(a)$. Section 3 is devoted to the proof Theorem \ref{th3} on the limiting behavior of nonnegative minimizers for $e(a)$ as $a\nearrow a^*$.

\section{Existence and Nonexistence of Minimizers}

This section is devoted to the proof of Theorems \ref{th1} and \ref{th2} with regard to the existence and nonexistence of minimizers for (\ref{e(a)}). We first recall the following compactness lemma \cite[Lemma 3.2]{C}.
\begin{lem}\label{lem2.1}
Assume $V(x)\in L^{\infty}_{loc}(\R^{N})$ and $\lim_{|x|\rightarrow\infty}V(x)=\infty$, then the embedding $\mathcal{H}\hookrightarrow L^{q}(\R^{N})$ is compact for $2\leq q <2_{s}^{*}$, where $2_s^*=\frac{2N}{N-2s}$ if $N>2s$ and $2_s^*=\infty$ if $N\leq2s$.
\end{lem}

We next discuss the decay estimate of $Q(x)$, where $Q(x)$ is a positive solution of the equation \eqref{Eq}.
\begin{lem}\label{lem2.3}
Under the assumptions of Theorem A, let $Q(x)>0$ be a solution of \eqref{Eq}. Then $Q(x)$ admits the following polynomially decay
\begin{equation}\label{ds}
|Q(x)|\leq C|x|^{-N-2s}\ \ \mbox{and} \ \ |\nabla Q(x)|\leq C|x|^{-N-1}\ \ \mbox{for}\ \ |x|> 1.
\end{equation}
\end{lem}
\vskip 0.1truein
\noindent\textbf{Proof.} Recall from \cite [Proposition 5]{RS} that
\begin{equation}\label{Z1}
|Q(x)|\leq C|x|^{-N},\ \ |\nabla Q(x)|\leq C|x|^{-N-1}\ \ \mbox{for}\ \ |x|>1\ \ \mbox{and}\ \ \  Q(x)\in C^{1}(\R^N\backslash \{0\}).
\end{equation}
In order to prove \eqref{ds}, it thus suffices to prove that 
\begin{equation*}\label{jZ1}
|Q(x)|\leq C|x|^{-N-2s}\ \ \mbox{for}\ \ |x|>1. 
\end{equation*}
Since $Q(x)>0$ satisfies \eqref{Eq}, we derive from \eqref{Z1} that there exists a large constant $R>0$ such that
\begin{equation}\label{Z2}
(-\Delta)^sQ+\frac{1}{2}Q\leq 0\ \ \mbox{in}\ \ B^c_{R}(0).
\end{equation}
For any $\lambda>0$, let $G_\lambda$ be the Fourier transform of $(|\xi|^{2s}+\lambda)^{-1}$.  Following \cite[Lemma C.1]{FLS}, one can obtain that there exists a positive $c>0$ such that $G_{\frac{1}{2}}(x)\geq c>0$ for $|x|\leq R$, and $G_{\frac{1}{2}}(x)\to 0$ as $|x|\to\infty$. Moreover, note from \cite [Proposition 4] {RS} that $Q(x)\in L^{\infty}(\R^N)$. Choosing $M:=\|Q\|_{L^{\infty}(\R^N)}c^{-1}$, we have
\begin{equation}\label{Z3}
Q(x)\leq \|Q\|_{L^{\infty}(\R^N)}=cM\leq M G_{\frac{1}{2}}\ \ \mbox{in}\ \ B_R(0).
\end{equation}
Define
\begin{equation}\label{Z4}
v(x):=M G_{\frac{1}{2}}(x)-Q(x),
\end{equation}
it then follows from \eqref{Z1} and \eqref{Z3} that 
$$v\geq 0\ \ \mbox{in}\ \ B_R(0)\backslash \{0\},\ \  v(x)\rightarrow 0\ \ \mbox{as}\ \ |x|\rightarrow\infty \ \ \mbox{and}\ \ v\in C(\R^N\backslash \{0\}).$$

We next prove that $v\geq0$ in $B_R^c(0)$ as well. On the contrary, suppose that there exists a $\bar x\in B_R^c(0)$ such that $v(\bar x)<0$. Noting that $v(x)\rightarrow 0$ as $|x|\rightarrow\infty$ and $v\geq 0$ in $B_R(0)$, we then deduce that there exists $x_0\in B_R^c(0)$ such that
\begin{equation}\label{z5A}
v(x_0)<0\ \ \mbox{and}\ \ v(x_0)<v(x)\ \ \mbox{for any}\ \ v\in C(\R^N\backslash \{0\}).
\end{equation}
 Using the singular integral expression for $(-\Delta)^{s}$, we deduce from \eqref{z5A} that
\begin{equation}\label{Z5}
((-\Delta)^{s}v)(x_0)<0\ \ \mbox{and}\ \ \ (-\Delta)^{s}v(x_0)+\frac{1}{2}v(x_0)<0.
\end{equation}
Since $G_{\frac{1}{2}}$ satisfies $(-\Delta)^{s}G_{\frac{1}{2}}+\frac{1}{2}G_{\frac{1}{2}}=\delta_0$ in $\R^N$, we derive from \eqref{Z2} that
\begin{equation}\label{Z6}
(-\Delta)^{s}v(x)+\frac{1}{2}v(x)\geq 0\ \ \mbox{in}\ \ B_R^c(0),
\end{equation}
which contradicts with \eqref{Z5}, we thus obtain that $v\geq 0$ in $B_R^c(0)$. Applying \cite[Lemma C.1]{FLS}, we derive from above that
\begin{equation}\label{Z7}
|Q(x)|\leq M G_{s,\frac{1}{2}}\leq \frac{C}{|x|^{N+2s}}\ \ \mbox{for any}\ \ x\in B_R^c(0).
\end{equation}
Together with the fact $Q(x)\in L^{\infty}(\R^N)$, it further implies that \eqref{ds} holds true. This completes the proof of Lemma \ref{lem2.3}.\qed

Due to the nonlocal nature of the operator $(-\Delta)^{s}$, in order to prove Theorems \ref{th1} and \ref{th2}, we next applying Lemma \ref{lem2.3} to estimate the Gagliardo (semi) norm of
\begin{equation}\label{c01}
Q_{\tau}(x)=\varphi(\tau^{-1}x)Q(x),\ \ \tau>0,
\end{equation}
where $Q(x)>0$ is a positive solution of \eqref{Eq}, and $\varphi\in C_{0}^{\infty}(\R^N, [0,1])$ is a nonnegative cut-off function  satisfying $\varphi(x)=1$ if $|x|\leq1$ and $\varphi(x)=0$ if $|x|>2$. 

From now on, for any $r>0$, we denote $B_r:=B_{r}(0)$ and $B_r^c:=\R^N\backslash B_{r}(0) $ for simplicity.

\begin{lem}\label{lem2.2}
Under the assumptions of Theorem \ref{th1}, we have 
\begin{equation*}\label{c2}
\iint_{\R^{2N}}\frac{|Q_{\tau}(x)-Q_{\tau}(y)|^2}{|x-y|^{N+2s}}dxdy\leq \iint_{\R^{2N}}\frac{|Q(x)-Q(y)|^2}{|x-y|^{N+2s}}dxdy+O(\tau^{-4s})\ \ \mbox{as}\ \ \tau\rightarrow\infty,
\end{equation*}
where $Q_\tau(x)> 0$ is defined by \eqref{c01}. 
\end{lem}

\vskip 0.05truein
\noindent\textbf{Proof.}
We first claim that the following assertions hold true:
\begin{enumerate}
\item For any $x\in\R^N$ and $y\in B^{c}_{\tau}$ with $|x-y|\leq \frac{\tau}{2}$, we have
\begin{equation}\label{c03}
|Q_{\tau}(x)-Q_{\tau}(y)|\leq C \tau^{-N-1}|x-y|;
\end{equation}
\item For any $x, y\in  B^{c}_{\tau}$,
\begin{equation}\label{c04}
|Q_{\tau}(x)-Q_{\tau}(y)|\leq C \tau^{-N-1}\min\{1,|x-y|\}
\end{equation}
for any $\tau\geq 2$ and some positive constant $C>0$.
\end{enumerate}

In fact, for any $x\in\R^N$ and $y\in B^{c}_{\tau}$ with $|x-y|\leq \frac{\tau}{2}$, let $
\xi$ be any point on the segment connecting $x$ and $y$, i.e.
\begin{equation*}
\xi=tx+(1-t)y\ \ \mbox{for some}\ \ t\in [0,1].
\end{equation*}
We then have
\begin{equation*}
|\xi|=|y+t(x-y)|\geq|y|-t|x-y|\geq\tau-\frac{t\tau}{2}\geq\frac{\tau}{2}\geq1.
\end{equation*}
Applying the differential mean value theorem, it follows from the polynomial decay of $\nabla Q$ in \eqref{ds} that
\begin{equation*}
|Q_{\tau}(x)-Q_{\tau}(y)|\leq C \tau^{-N-1}|x-y|,
\end{equation*}
and the assertion (1) is proved. We next focus on proving assertion (2), for any $x, y\in  B^{c}_{\tau}$, note that if $|x-y|\leq 1$, then the assertion (2) follows from (1) because of $\tau\geq 2$. For $|x-y|>1$, we derive from \eqref{ds} and \eqref{c01} that
\begin{equation*}
|Q_{\tau}(x)-Q_{\tau}(y)|\leq|Q(x)|+|Q(y)|\leq C \tau^{-N-2s}\leq C \tau^{-N-1}.
\end{equation*}
As a result, we conclude that $|Q_{\tau}(x)-Q_{\tau}(y)|\leq C \tau^{-N-1}\min\{1,|x-y|\}$ and the claim is therefore complete.

Now we denote
\begin{equation*}
\mathbb{E}:=\Big\{(x,y)\in\R^{N}\times\R^N: x\in B_{\tau},\ \ y\in B^{c}_{\tau} \ \ \mbox{and}\ \ |x-y|\leq \frac{\tau}{2}\Big\},
\end{equation*}
and
\begin{equation*}
\mathbb{F}:=\Big\{(x,y)\in\R^{N}\times\R^N: x\in B_{\tau},\ \ y\in B^{c}_{\tau} \ \ \mbox{and}\ \ |x-y|> \frac{\tau}{2}\Big\}.
\end{equation*}
It then follows from \eqref{c01} that
\begin{align}\label{c5}
I:=& \iint_{\R^{2N}}\frac{|Q_{\tau}(x)-Q_{\tau}(y)|^2}{|x-y|^{N+2s}}dxdy\nonumber\\
=&\iint_{B_{\tau}\times B_{\tau}}\frac{|Q(x)-Q(y)|^2}{|x-y|^{N+2s}}dxdy+2\iint_{\mathbb{E}}\frac{|Q_{\tau}(x)-Q_{\tau}(y)|^2}{|x-y|^{N+2s}}dxdy\\
&+2\iint_{\mathbb{F}}\frac{|Q_{\tau}(x)-Q_{\tau}(y)|^2}{|x-y|^{N+2s}}dxdy+\iint_{B^{c}_{\tau}\times B^{c}_{\tau}}\frac{|Q_{\tau}(x)-Q_{\tau}(y)|^2}{|x-y|^{N+2s}}dxdy\nonumber\\
:=&I_1+2I_2+2I_3+I_4.\nonumber
\end{align}
Next, we estimate the right hand terms of \eqref{c5}. For $I_2$, we deduce from \eqref{c03} that
\begin{equation}\label{c6}
\begin{split}
I_2&=\iint_{\mathbb{E}}\frac{|Q_{\tau}(x)-Q_{\tau}(y)|^2}{|x-y|^{N+2s}}dxdy\\
&\leq C\tau^{-2N-2}\iint_{\mathbb{E}}\frac{|x-y|^2}{|x-y|^{N+2s}}dxdy\\
&\leq C\tau^{-2N-2}\int_{|x|\leq\tau}dx\int_{|\xi|\leq\frac{\tau}{2}}\frac{1}{|\xi|^{N+2s-2}}d\xi=O(\tau^{-N-2s})\ \ \mbox{as}\ \ \tau\rightarrow\infty.
\end{split}
\end{equation}

As for $I_3$, we obtain from \eqref{c01} that $Q_{\tau}(x)=Q(x)$ for any
$x\in B_{\tau}$. For any $(x,y)\in \mathbb{F}$, it then holds that
\begin{equation*}\label{c7}
\begin{split}
\big|Q_{\tau}(x)-Q_{\tau}(y)\big|^2=&\big|(Q(x)-Q(y))+(Q(y)-Q_{\tau}(y))\big|^2\\
\leq&\big|Q(x)-Q(y)|^2+|Q(y)-Q_{\tau}(y)\big|^2\\
&+2\big|(Q(x)-Q(y))(Q(y)-Q_{\tau}(y))\big|,
\end{split}
\end{equation*}
which further implies that
\begin{align}\label{c8}
I_3=&\iint_{\mathbb{F}}\frac{|Q_{\tau}(x)-Q_{\tau}(y)|^2}{|x-y|^{N+2s}}dxdy\nonumber\\
\leq&\iint_{\mathbb{F}}\frac{|Q(x)-Q(y)|^2}{|x-y|^{N+2s}}dxdy+\iint_{\mathbb{F}}\frac{|Q(y)-Q_{\tau}(y)|^2}{|x-y|^{N+2s}}dxdy\\
&+2\iint_{\mathbb{F}}\frac{|Q(x)-Q(y)||Q(y)-Q_{\tau}(y)|}{|x-y|^{N+2s}}dxdy:=I_{31}+I_{32}+2I_{33}.\nonumber
\end{align}
By the polynomial decay \eqref{ds} of $Q$, it yields from \eqref{c01} that
\begin{equation}\label{c9}
\begin{split}
I_{32}&=\iint_{\mathbb{F}}\frac{|Q(y)-Q_{\tau}(y)|^2}{|x-y|^{N+2s}}dxdy\\
&\leq4\iint_{\mathbb{F}}\frac{|Q(y)|^2}{|x-y|^{N+2s}}dxdy\\
&\leq C\tau^{-2N-4s}\iint_{\mathbb{F}}\frac{1}{|x-y|^{N+2s}}dxdy\\
&\leq C\tau^{-2N-4s}\int_{|x|\leq\tau}dx\int_{|\xi|>\frac{\tau}{2}}\frac{1}{|\xi|^{N+2s}}d\xi=O(\tau^{-N-6s})\ \ \mbox{as}\ \ \tau\rightarrow\infty.
\end{split}
\end{equation}

For $I_{33}$, recall from \cite[Proposition 4]{RS} that $Q(x)\in L^{\infty}(\R^N)$. We then deduce from \eqref{ds} that
\begin{equation}\label{c9A}
|Q(x)Q(y)|\leq C\tau^{-N-2s}\ \ \mbox{for any}\ \ (x,y)\in \mathbb{F}.
\end{equation}
It then follows from \eqref{c01} and \eqref{c9A} that
\begin{equation}\label{c010}
\begin{split}
\iint_{\mathbb{F}}\frac{|Q(y)-Q_{\tau}(y)||Q(x)|}{|x-y|^{N+2s}}dxdy&\leq 2\iint_{\mathbb{F}}\frac{|Q(y)||Q(x)|}{|x-y|^{N+2s}}dxdy\\
&\leq C \tau^{-N-2s}\iint_{\mathbb{F}}\frac{1}{|x-y|^{N+2s}}dxdy\\
&\leq C \tau^{-N-2s}\int_{|x|\leq\tau}dx\int_{|\xi|>\frac{\tau}{2}}\frac{1}{|\xi|^{N+2s}}d\xi\\
&=O(\tau^{-4s})\ \ \mbox{as}\ \ \tau\rightarrow\infty.
\end{split}
\end{equation}
Moreover, we derive from \eqref{ds} and \eqref{c01}  that
\begin{equation}\label{c11}
\begin{split}
\iint_{\mathbb{F}}\frac{|Q(y)-Q_{\tau}(y)||Q(y)|}{|x-y|^{N+2s}}dxdy&\leq 2\iint_{\mathbb{F}}\frac{|Q(y)|^2}{|x-y|^{N+2s}}dxdy\\
&\leq C \tau^{-2N-4s} \iint_{\mathbb{F}}\frac{1}{|x-y|^{N+2s}}dxdy\\
&\leq C \tau^{-2N-4s}\int_{|x|\leq\tau}dx\int_{|\xi|>\frac{\tau}{2}}\frac{1}{|\xi|^{N+2s}}d\xi\\
&=O(\tau^{-N-6s})\ \ \mbox{as}\ \ \tau\rightarrow\infty.
\end{split}
\end{equation}
Together with \eqref{c010}, it further yields that
\begin{equation}\label{c12}
\begin{split}
I_{33}&\leq \iint_{\mathbb{F}}\frac{|Q(y)-Q_{\tau}(y)||Q(x)|}{|x-y|^{N+2s}}dxdy+\iint_{\mathbb{F}}\frac{|Q(y)-Q_{\tau}(y)||Q(y)|}{|x-y|^{N+2s}}dxdy\\
&\leq O(\tau^{-4s})\ \ \mbox{as}\ \ \tau\rightarrow\infty.
\end{split}
\end{equation}
Consequently, we deduce that
\begin{equation*}\label{jc13}
I_3\leq \iint_{\mathbb{F}}\frac{|Q(x)-Q(y)|^2}{|x-y|^{N+2s}}dxdy+O(\tau^{-4s})\ \ \mbox{as}\ \ \tau\rightarrow\infty.
\end{equation*}

Finally, we estimate $I_4$ as follows. We deduce from \eqref{c01} and \eqref{c04} that
\begin{equation}\label{c14}
\begin{split}
I_4&=\iint_{B^{c}_{\tau}\times B^{c}_{\tau}}\frac{|Q_{\tau}(x)-Q_{\tau}(y)|^2}{|x-y|^{N+2s}}dxdy\\
&\leq C \tau^{-2N-2}\iint_{B_{2\tau}\times \R^N}\frac{\min\{1,|x-y|^2\}}{|x-y|^{N+2s}}dxdy\\
&=C \tau^{-2N-2}\int_{|x|\leq2\tau}dx\Big(\int_{|\xi|\leq 1}\frac{1}{|\xi|^{N+2s-2}}d\xi+\int_{|\xi|\geq 1}\frac{1}{|\xi|^{N+2s}}d\xi\Big)\\
&=O(\tau^{-N-2})\ \ \mbox{as}\ \  \tau\rightarrow\infty.
\end{split}
\end{equation}
In view of the above facts, we conclude that
\begin{align}\label{c15}
I&:= \iint_{\R^{2N}}\frac{|Q_{\tau}(x)-Q_{\tau}(y)|^2}{|x-y|^{N+2s}}dxdy\nonumber\\
&\leq\iint_{B_{\tau}\times B_{\tau}}\frac{|Q(x)-Q(y)|^2}{|x-y|^{N+2s}}dxdy+2\iint_{\mathbb{F}}\frac{|Q(x)-Q(y)|^2}{|x-y|^{N+2s}}dxdy+O(\tau^{-4s})\\
&\leq\iint_{\R^{2N}}\frac{|Q(x)-Q(y)|^2}{|x-y|^{N+2s}}dxdy+O(\tau^{-4s}) \ \ \mbox{as}\ \ \tau\rightarrow\infty.\nonumber
\end{align}
The proof of Lemma \ref{lem2.2} is thus complete.\qed

In order to complete the proof of Theorem \ref{th1}, we shall also give the following lemma.
\begin{lem}\label{A.111}
Suppose $\Omega$ is a bounded domain in $\R^N$ or $\Omega=\R^N$, and assume that $u_{n}\rightarrow u$ strongly in $L^{q}(\Omega)$ for any $q\in[2,2_s^*)$. Let $s\in (\frac{1}{2},1)$, $N>2s$, $0<b<\min\{\frac{N}{2},1\}$, $\beta=\sqrt{\frac{2s-b}{N}}$ and $f_n(x)=|u_{n}(x)|^{2+2\beta^2}-|u(x)|^{2+2\beta^2}$, then $f_n\to 0$ in $L^m(\Omega)$ for any  $m\in\big[1,\frac{N^2}{(N-2s)(N+2s-b)}\big)$ if $N>2s$; $f_n\to 0$ in $L^m(\Omega)$ for any $m\in [1,\infty)$ if $N\leq 2s$.
\end{lem}
Since the proof of Lemma \ref{A.111} is long and similar to \cite[Lemma B.1 ]{LZ}, we omit it for simplicity.

Based on above Lemmas, applying the Gagliardo-Nirenberg inequality (\ref{GN}), we shall complete the proofs of Theorems \ref{th1} and \ref{th2}.

\vskip 0.05truein
\noindent{\bf Proof of Theorem \ref{th1}.}
(1). Firstly, we prove that $e(a)$ has at least one minimizer for any $0< a< a^{\ast}=\|Q\|^{2\beta^2}_{2}$.
Indeed, for any $0< a<a^{\ast}$ and $u\in \mathcal{H}$ with $\|u\|_{L^2(\R^N)}=1$, we derive from the assumption $(V_1)$ and the Gagliardo-Nirenberg inequality (\ref{GN}) that
\begin{align}\label{Eu}
E_{a}(u)&=\int_{\R^{N}}(|(-\Delta)^{\frac{s}{2}}u(x)|^{2}+V(x)|u(x)|^{2})dx-\frac{a}{1+\beta^2}\int_{\R^{N}}|u(x)|^{2+2\beta^2}|x|^{-b}dx  \nonumber\\
&\geq(1-\frac{a}{a^*})\int_{\R^{N}}|(-\Delta)^{\frac{s}{2}} u(x)|^{2}dx+\int_{\R^{N}}V(x)|u(x)|^{2}dx \\
&\geq(1-\frac{a}{a^*})\int_{\R^{N}}|(-\Delta)^{\frac{s}{2}} u(x)|^{2}dx\geq 0,\nonumber
\end{align}
which further implies that $e(a)\geq0$ for $0<a<a^*$. Then one can choose a minimizing sequence $\{u_{n}\}\!\subset\!\mathcal{H}$ satisfying $\|u_{n}\|_{L^2(\R^N)}\!=\!1$ and $\lim\limits_{n\rightarrow\infty}\!\!E_{a}(u_{n})\!\!=\!\!e(a)$. Due to (\ref{Eu}), we obtain that \!$\int_{\R^{N}}\!|(-\Delta)^{\frac{s}{2}} u_{n}(x)|^{2}dx$ and \!$\int_{\R^{N}}\!\!V(x)|u_{n}|^{2}dx$ are bounded uniformly for $n$. Following Lemma \ref{lem2.1}, passing  to a subsequence if necessary, we may assume that
$$u_{n}\rightharpoonup u\,\ \hbox{weakly in} \,\ \mathcal{H},\,\  u_{n}\rightarrow u\,\ \hbox{strongly in}\,\ L^{q}(\R^{N}),\,\ 2\leq q<2_s^* $$
for some $u\in \mathcal{H}$.
Thus we deduce that $\int_{\R^{N}}|u(x)|^{2}dx=1$.

On the other hand, for any $t\in (\frac{N}{N-b},\infty)$ and $R>0$, note that any function $f\in L^1(\R^N)\cap L^t(\R^N)$ satisfies
\begin{align}\label{02.4}
\inte f|x|^{-b}dx
\leq& \|f\|_{L^{t}(B_R(0))}\big\||x|^{-b}\big\|_{L^{\bar t}(B_R(0))}+\frac{1}{R^b}\|f\|_{L^{1}(\R^N)} \nonumber\\
\leq& C(R)\|f\|_{L^{t}(B_R(0))}+\frac{1}{R^b}\|f\|_{L^{1}(\R^N)},
\end{align}
where $\bar t\!=\!\frac{t}{t-1}$ satisfies $1\!<\!\bar t\!<\!\frac{N}{b}$. Therefore, taking $f_n\!=\!|u_{n}(x)|^{2+2\beta^2}\!-\!|u(x)|^{2+2\beta^2}$, we deduce from Lemma \ref{A.111} that $f_n\to 0$ in $L^1((B_R(0)))\cap L^t((B_R(0)))$ for $t\in\big[1,\frac{N^2}{(N-2s)(N+2s-b)}\big)$ if $N>2s$ and $t\in [1,\infty)$ if $N\leq 2s$. Let $R\rightarrow\infty$ in \eqref{02.4}, it then follows that
\begin{equation}\label{2.1}
\lim_{n\rightarrow\infty}\int_{\R^{N}}|x|^{-b}|u_{n}(x)|^{2+2\beta^2}dx=\int_{\R^{N}}|x|^{-b}|u(x)|^{2+2\beta^2}dx.
\end{equation}
Consequently, by the weak lower semicontinuity, we have
\begin{equation*}
e(a)=\lim_{n\to\infty}E_{a}(u_n)\\ \geq E_{a}(u)\geq e(a).
\end{equation*}
Therefore, $u$ must be a minimizer of $e(a)$ and Theorem \ref{th1} (1) is proved.

(2). We next prove that there is no minimizer for $e(a)$ as $a>a^{\ast}$. For any $\tau>0$, set
\begin{equation}\label{st}
u_{\tau}(x):=A_{\tau}\frac{\tau^{\frac{N}{2}}}{\|Q\|_{2}} Q_{\tau}(\tau x),
\end{equation}
where $A_{\tau}>0$ is chosen so that $\int_{\R^{N}}|u_{\tau}(x)|^{2}dx=1$ and $Q_{\tau}(x)$ is defined by \eqref{c01}. By scaling, we deduce from the polynomial decay of $Q$ in (\ref{ds}) that
\begin{equation}\label{A}
\frac{1}{A^{2}_{\tau}}=\frac{1}{\|Q\|^{2}_{2}}\int_{\R^{N}}Q^{2}(x)\varphi^{2}(\tau^{-1}x)=1+O(\tau^{-N-4s}) \ \ \mbox{as}\ \ \tau\rightarrow\infty.
\end{equation}
According to (\ref{DS}), (\ref{ds}) and Lemma \ref{lem2.2}, we also obtain that
\begin{align}\label{Gj}
&\int_{\R^{N}}|(-\Delta)^{\frac{s}{2}} u_{\tau}(x)|^{2}dx-\frac{a}{1+\beta^2}\int_{\R^{N}}|x|^{-b}u_{\tau}^{2+2\beta^2}dx
\nonumber\\=&\frac{\tau^{2s}A^{2}_{\tau}}{\|Q\|^{2}_{2}}\int_{\R^{N}}|(-\Delta)^{\frac{s}{2}} Q_{\tau}|^{2}dx-\frac{a}{1+\beta^2}\int_{\R^{N}}|x|^{-b}u_{\tau}^{2+2\beta^2}dx
\nonumber\\ \leq &\frac{\tau^{2s}A^{2}_{\tau}}{\|Q\|^{2}_{2}}\int_{\R^{N}}|(-\Delta)^{\frac{s}{2}} Q|^{2}dx-\frac{a}{1+\beta^2}\frac{\tau^{2s}A^{2+2\beta^2}_{\tau}}{\|Q\|^{2+2\beta^2}_{2}}\int_{\R^{N}}|x|^{-b}Q^{2+2\beta^2}dx+O(\tau^{-2s})
\\=&\frac{\tau^{2s}A^{2}_{\tau}}{\|Q\|^{2}_{2}}\int_{\R^{N}}|(-\Delta)^{\frac{s}{2}} Q|^{2}dx-\frac{a\tau^{2s}A^{2+2\beta^2}_{\tau}}{\|Q\|^{2+2\beta^2}_{2}}\int_{\R^{N}}|(-\Delta)^{\frac{s}{2}} Q|^{2}dx+O(\tau^{-2s})
\nonumber\\=&\frac{\tau^{2s}}{\|Q\|^{2}_{2}}\big(1-\frac{a}{a^*}\big)\int_{\R^{N}}|(-\Delta)^{\frac{s}{2}} Q|^{2}dx+O(\tau^{-2s})\ \ \mbox{as}\ \ \tau\rightarrow\infty.	\nonumber
\end{align}

On the other hand, since the function $V(x)\varphi^{2}(x)$ is bounded and has compact support, it follows from \cite[Theorem 1.8]{LE} that
\begin{equation}\label{jx}
\lim\limits_{\tau\rightarrow\infty}\int_{\R^{N}}V(x)|u_{\tau}(x)|^{2}dx=V(0).
\end{equation}
In view of \eqref{Gj} and \eqref{jx}, we then derive that
\begin{align}\label{zh}
e(a)\leq &E_{a}(u_{\tau}(x)) \nonumber\\
=&\int_{\R^{N}}|(-\Delta)^{\frac{s}{2}} u_{\tau}(x)|^{2}dx-\frac{a}{1+\beta^2}\int_{\R^{N}}|x|^{-b}u^{2+2\beta^2}_{\tau}(x)dx\nonumber\\
&+\int_{\R^{N}}V(x)|u_{\tau}(x)|^{2}dx\\
=&\frac{\tau^{2s}}{\|Q\|^{2}_{2}}(1-\frac{a}{a^*})\int_{\R^{N}}|(-\Delta)^{\frac{s}{2}} Q|^{2}dx+\int_{\R^{N}}V(x)|u_{\tau}(x)|^{2}dx+O(\tau^{-2s})\nonumber\\
\rightarrow&-\infty\ \ \hbox{as}\ \  \tau\to\infty.\nonumber
\end{align}
Therefore, if $a>a^{\ast}$, $e(a)=-\infty$, which implies the nonexistence of minimizers as $a>a^{\ast}$. Finally, applying Theorem \ref{th1} (1) and the Gagliardo-Nirenberg inequality \eqref{GN}, it is obviously that $e(a)\!>\!0$ for $0\!<\!a\!<\!a^*$. This completes the proof of Theorem \ref{th1}.\qed

\vskip 0.05truein
\noindent{\bf Proof of Theorem \ref{th2}.} (1). Applying the Gagliardo-Nirenberg inequality \eqref{GN}, one can check that $e(a^*)\geq 0$ is bounded from below, and thus there exists a minimizing sequence $\{u_{n}\}$ of $e(a^{\ast})$ such that $e(a^*)=\lim_{n\to\infty}E_a(u_n)$.
To prove Theorem \ref{th2} (1), similar to the proof of Theorem \ref{th1}, it is enough to prove that $\{u_{n}\}$ is bounded uniformly in $\mathcal{H}$. On the contrary, suppose that $\{u_{n}\}$ is unbounded in $\mathcal{H}$, then there exists a subsequence, still denoted by $\{u_{n}\}$, of $\{u_{n}\}$ such that $\|u_{n}\|_{\mathcal{H}}\rightarrow\infty$ as $n\rightarrow\infty$. Moreover, by the Gagliardo-Nirenberg inequality \eqref{GN}, we have
\begin{equation*}\label{V}
\int_{\R^{N}}V(x)|u_{n}(x)|^{2}dx\leq E_{a^{\ast}}(u_{n})\leq e(a^{\ast})+1\ \ \mbox{for\ \ sufficiently\ \ large}\ \ n>0.
\end{equation*}
Noting that $\{u_{n}\}$ is unbounded in $\mathcal{H}$, taking a subsequence if necessary, we have
\begin{equation}\label{U}
\int_{\R^{N}}|(-\Delta)^{\frac{s}{2}} u_{n}(x)|^{2}dx\rightarrow\infty\,\ \hbox{as}\,\ n\to\infty.
\end{equation}
Define
\begin{equation*}\label{wq}
\varepsilon^{-2s}_{n}:=\int_{\R^{N}}|(-\Delta)^{\frac{s}{2}}  u_{n}(x)|^{2}dx,
\end{equation*}
it then follows from (\ref{U}) that $\varepsilon_{n}\rightarrow 0$ as $n\rightarrow\infty$. We next define
\begin{equation*}\label{omg}
w_{n}(x):=\varepsilon^{\frac{N}{2}}_{n}u_{n}(\varepsilon_{n}x),
\end{equation*}
so that
\begin{equation}\label{2gj}
\int_{\R^{N}}|(-\Delta)^{\frac{s}{2}}  w_{n}(x)|^{2}dx=\int_{\R^{N}}| w_{n}(x)|^{2}dx=1.
\end{equation}
Therefore $\{w_n\}$ is bounded uniformly for $n$ in $H^s(\R^N)$. By Sobolev embedding, passing to a subsequence if necessary, there exists a function $\bar{w}_0\in H^s(\R^3)$ such that
\begin{equation}\label{s1}
w_{n}\rightharpoonup \bar{w}_0 \,\ \hbox{weakly in} \,\ H^s(\R^N),\,\  w_{n}\rightarrow \bar{w}_0\,\ \hbox{strongly in}\,\ L_{loc}^{q}(\R^{N}),\,\ 2\leq q<2_s^*.
\end{equation}

We next claim that
\begin{equation}\label{2.2}
\lim_{n\to\infty}w_n=\bar{w}_0\ \ \hbox{strongly in $H^s(\R^N)$},
\end{equation}
where $\bar{w}_0$ is an optimizer of  the Gagliardo-Nirenberg inequality \eqref{GN}. Actually, by the definition of $w_n$, we obtain from (\ref{GN}) that
\begin{align}\label{wgx}
C\geq &e(a^*)=\lim_{n\rightarrow\infty}E_{a^*}(u_{n})
\nonumber\\=&\lim_{n\rightarrow\infty}\Big\{\int_{\R^{N}}\Big(|(-\Delta)^{\frac{s}{2}} u_{n}|^{2}-\frac{a^*}{1+\beta^2}|x|^{-b}u_{n}^{2+2\beta^2}+V(x)|u_{n}(x)|^{2}\Big)dx\Big\}
\nonumber\\
=&\lim_{n\rightarrow\infty}\Big\{\frac{1}{\varepsilon^{2s}_{n}}\int_{\R^{N}}\Big[|(-\Delta)^{\frac{s}{2}} w_{n}|^{2}-\frac{a^{\ast}}{1+\beta^2}|x|^{-b}w_{n}^{2+2\beta^2}\Big]dx\\
&+\int_{\R^{N}}V(\varepsilon_{n}x)|w_{n}(x)|^{2}dx\Big\}
\geq 0.\nonumber
\end{align}
Therefore, we derive from \eqref{2gj}, the fact $\eps_n\to 0$ as $n\to\infty$ and above that
\begin{equation*}\label{wnp}
\lim_{n\to\infty}\frac{a^{\ast}}{1+\beta^2}\int_{\R^{N}}|x|^{-b}w_{n}^{2+2\beta^2}(x)dx=\lim_{n\to\infty}\int_{\R^{N}}|(-\Delta)^{\frac{s}{2}}  w_{n}(x)|^{2}dx=1.
\end{equation*}
On the other hand, applying Lemma \ref{A.111} and \eqref{02.4}, let $R\rightarrow\infty$ in \eqref{02.4} we deduce that
\begin{equation}\label{eq2.20}
\frac{a^{\ast}}{1+\beta^2}\int_{\R^{N}}|x|^{-b}\bar{w}_0^{2+2\beta^2}(x)dx=\lim_{n\to\infty}\frac{a^{\ast}}{1+\beta^2}\int_{\R^{N}}|x|^{-b}w_{n}^{2+2\beta^2}(x)dx=1,
\end{equation}
which further implies $\bar{w}_0\not\equiv0$. We then deduce from (\ref{GN}) and \eqref{s1} that
\begin{align*}\label{wgx1}
0&= \lim_{n\rightarrow\infty}\varepsilon_n^{2s}e(a^*)
\nonumber\\
&\geq\lim_{n\rightarrow\infty}\int_{\R^{N}}\Big[|(-\Delta)^{\frac{s}{2}} w_{n}|^{2}-\frac{a^{\ast}}{1+\beta^2}|x|^{-b}w_{n}^{2+2\beta^2}\Big]dx
\nonumber\\&\geq\int_{\R^{N}}|(-\Delta)^{\frac{s}{2}} \bar{w}_0(x)|^{2}dx-\frac{a^{\ast}}{1+\beta^2}\int_{\R^{N}}|x|^{-b}\bar{w}_0^{2+2\beta^2}(x)dx
\nonumber\\&\geq\Big[1-\Big(\int_{\R^{N}}
\bar{w}^2_0(x)dx\Big)^{\beta^2}\Big]\int_{\R^{N}}|(-\Delta)^{\frac{s}{2}} \bar w_{0}(x)|^{2}dx.
\end{align*}
Since $\|\bar{w}_0\|_{L^2(\R^N)}\leq 1$, the above inequality further implies that
\begin{equation}\label{2.3}
\int_{\R^{N}}|\bar{w}_0(x)|^{2}dx\!=\!1,\,\ \int_{\R^{N}}|(-\Delta)^{\frac{s}{2}} \bar{w}_0(x)|^{2}dx\!=\!\frac{a^{\ast}}{1+\beta^2}\int_{\R^{N}}|x|^{-b}\bar{w}_0^{2+2\beta^2}(x)dx\!=\!1.
\end{equation}

In view of above facts, by the norm preservation, we conclude that $w_{n}$ converges to $\bar{w}_0$ strongly in $H^{s}(\R^{N})$ as $n\to\infty$, where $\bar{w}_0$ is an optimizer of  the Gagliardo-Nirenberg inequality \eqref{GN} and \eqref{2.2} is thus proved. Applying Fatou's lemma,  we then derive from \eqref{GN} and \eqref{2.2} that
\begin{equation*}\label{md1}
e(a^{\ast})\geq\lim_{n\rightarrow\infty}\int_{\R^{N}}V(\varepsilon_{n} x)|w_{n}(x)|^{2}dx\geq \int_{\R^{N}}V(0)|\bar{w}_0(x)|^{2}dx=V(0),
\end{equation*}
which contradicts with the assumption  $e(a^{\ast})<V(0)$. Hence the minimizing sequence $\{u_{n}\}$ of $e(a^*)$ is bounded uniformly in $\mathcal{H}$ and Theorem \ref{th2} (1) is thus proved.

(2). We next consider the case $a=a^{\ast}$ and $V(0)=0$. Note that $e(a^*)\geq 0$. On the other hand, one can use the same test function as that in (\ref{st}) to derive $e(a^{\ast})\leq V(0)=0$, and thus $e(a^{\ast})=0$. Now suppose that there exists a minimizer $u\in\mathcal{H}$ satisfying $\inte |u|^2dx=1$ at $a=a^{\ast}$. Then we have
\begin{equation*}\label{1}
\int_{\R^{N}}V(x)|u(x)|^{2}dx=0=\inf\limits_{x\in\R^{N}}V(x),
\end{equation*}
as well as
\begin{equation*}\label{2}
\int_{\R^{N}}|(-\Delta)^{\frac{s}{2}}u(x)|^{2}dx=\frac{a^{\ast}}{1+\beta^2}\int_{\R^{N}}|x|^{-b}|u(x)|^{2+2\beta^2}dx.
\end{equation*}
This is a contradiction, since the first equality implies that $u$ has compact support. However, we deduce from the second equality that $u>0$ is an optimizer of $\eqref{GN}$ and solves the equation \eqref{Eq}. Moreover, taking the same test function as in (\ref{st}), we have $V(0)=\inf_{x\in\R^{N}}V(x)\leq\lim_{a\nearrow a^{\ast}}e(a)\leq V(0)$ and hence $\lim_{a\nearrow a^{\ast}}e(a)=V(0)=e(a^*)$. This completes the proof of Theorem \ref{th2}.\qed

\section{Limiting Behavior of Minimizers as $a\nearrow a^{\ast}$}
In this section, we shall prove Theorem \ref{th3} on the limiting  behavior of nonnegative minimizers of $e(a)$ as $a\nearrow a^{\ast}$. Let $u_{a}$ be a nonnegative minimizer of $e(a)$, by variational theory,  $u_{a}$ satisfies the following Euler-Lagrange equation
\begin{equation}\label{EL}
(-\Delta)^{s} u_{a}(x)+V(x)u_{a}(x)=\mu_a u_{a}(x)+a u_{a}^{1+2\beta^2}(x)|x|^{-b}\ \ \mbox{in}\ \ \R^{N},
\end{equation}
where $\mu_{a}\in\R$ is a suitable Lagrange multiplier satisfying
\begin{equation}\label{2.8}
\mu_a=e(a)-\frac{a\beta^2}{1+\beta^2}\inte |x|^{-b}|u_a|^{2+2\beta^2}dx.
\end{equation}
Define
\begin{equation}\label{yp}
\varepsilon_{a}:=\Big(\int_{\R^{N}}|(-\Delta)^{\frac{s}{2}}u_{a}|^{2}dx\Big)^{-\frac{1}{2s}},
\end{equation}
and
\begin{equation}\label{xpa}
w_{a}(x):=\varepsilon^{\frac{N}{2}}_{a}u_{a}(\varepsilon_{a}x).
\end{equation}
Using the simialr argument of \cite[Lemma B.2 ]{LZ},  one can check that $w_a$ satisfies the following lemma.
\begin{lem}\label{A.1}
Let $s\in(\frac{1}{2},1)$, $N>2s$, $0<b<\min\{\frac{N}{2},1\}$ and $\beta^2=\frac{2s-b}{N}$. Suppose $w_a(x)\in H^s(\R^N)$, then $w^q_a(x)|x|^{-b}\in L^m(\R^N)$, where $m\in\big(\frac{2N}{Nq+2b}, \frac{2_s^*N}{Nq+2_s^*b}\big)$.
\end{lem}

Based on above facts, we then complete the proof of Theorem \ref{th3}.
\vskip 0.05truein
\noindent\textbf{Proof of Theorem \ref{th3}.} We first claim that $\varepsilon_{a}\rightarrow 0$ as  $a\nearrow a^{\ast}$. Actually, suppose on the contrary, there exists a sequence $\{a_{k}\}$ satisfying $a_{k}\nearrow a^{\ast}$ as $k\rightarrow\infty$, such that the sequence $\{u_{a_{k}}\}$ is bounded uniformly for $k$ in $H^s(\R^N)$. By the Gagliardo-Nirenberg inequality \eqref{GN}, we obtain that $\inte V(x)|u_{a_k}|^2dx\leq e(a_k)\leq C$ uniformly for $k$ and hence $\{u_{a_{k}}\}$ is bounded uniformly in $\mathcal{H}$.
Following Lemma \ref{lem2.1},  there exist a subsequence, still denoted by $\{a_{k}\}$, of $\{a_{k}\}$ and $u_{0}\in\mathcal{H}$ such that
$${u_{a_{k}}}\rightharpoonup u_{0}\,\ \hbox{weakly in}\,\ \mathcal{H}\,\ \hbox{and}\,\   {u_{a_{k}}}\rightarrow u_{0}\,\ \hbox{strongly in}\,\  L^{q}(\R^{N})\,(2\leq q <2_s^*)\ \ \mbox{as}\ \ k\to\infty.$$
Combining the above convergence with \eqref{2.1} further gives that
\begin{equation*}
0=e(a^{\ast})\leq E_{a^{\ast}}(u_{0})\leq \lim\limits_{k\rightarrow\infty}E_{a_{k}}(u_{a_{k}})=\lim\limits_{k\rightarrow\infty}e(a_{k})=0,
\end{equation*}
and thus $u_{0}$ is a minimizer of $e(a^{\ast})$, which however contradicts with Theorem \ref{th2} (2). Therefore we conclude that $\varepsilon_{a}\rightarrow 0$ as $a\nearrow a^{\ast}$.

We then prove that
\begin{equation}\label{2.5}
w_{a}\rightarrow w_0\,\ \hbox{strongly in}\,\ H^s(\R^N)\,\ \hbox{as}\,\ a\nearrow a^*,
\end{equation}
where $w_0$ is an optimizer of  the Gagliardo-Nirenberg inequality \eqref{GN}. Noting from Theorem \ref{th2} (2) that $0=e(a^*)=\lim_{a\nearrow a^*}e(a)$, thus $\{u_a\}$ is also a minimizing sequence of $e(a^*)$. Moreover, we obtain from \eqref{yp} and \eqref{xpa} that
\[
\inte |(-\Delta)^{\frac{s}{2}} w_a|^2dx=\inte |w_a|^2dx=1.
\]
Hence $\{w_a\}$ is bounded  in $H^s(\R^N)$ uniformly for $a$, and passing to a subsequence if necessary, there exists a function $w_0\in H^s(\R^N)$ such that
$$w_{a}\rightharpoonup w_0\,\ \hbox{weakly in} \,\ H^s(\R^N),\,\  w_{a}\rightarrow w_0\,\ \hbox{strongly in}\,\ L_{loc}^{q}(\R^{N}),\,\ 2\leq q<2_s^*. $$
Note that the proof of the claim $\eqref{2.2}$ does not rely on the assumption $0\leq e(a^*)< V(0)$. Thus, the same argument of proving \eqref{2.2} then yields that \eqref{2.5} holds true.

Next, we prove \eqref{jbp} holds true. Recall from \eqref{EL} and \eqref{xpa} that $w_a$ satisfies the following equation
\begin{equation}\label{2.9}
(-\Delta )^sw_a+\eps_a^{2s}V(\eps_ax)w_a=\mu_a\eps_a^{2s}w_a+a w_a^{1+2\beta^2} |x|^{-b}\ \ \hbox{in}\,\ \R^N,
\end{equation}
where $\mu_a\in\R$ is a suitable Lagrange multiplier  satisfying \eqref{2.8}. We derive from \eqref{2.3} and \eqref{2.8} that
\begin{align}\label{2.10}
\eps_a^{2s}\mu_a&=\eps_a^{2s}e(a)-\eps_a^{2s}\frac{a\beta^2}{1+\beta^2}\inte |x|^{-b}|u_a|^{2+2\beta^2}dx
\nonumber\\
&=\eps_a^{2s}e(a)-\frac{a\beta^2}{1+\beta^2}\inte |x|^{-b}|w_a|^{2+2\beta^2}dx\to -\beta^2\,\ \hbox{as}\,\ a\nearrow a^*.
\end{align}
Hence for $a\nearrow a^*$, it follows from \eqref{2.9} that
\begin{equation}\label{R1}
(-\Delta)^s w_{a}-c(x)w_{a}(x)\leq0\ \ \mbox{in}\ \ \R^N,\ \ \mbox{where}\ \ c(x)=aw^{2\beta^2}_{a}|x|^{-b}.
\end{equation}
In view of above facts, we shall prove \eqref{jbp} through the following two steps.

\vskip 0.05truein

\noindent{\em  Step 1.} We prove that there exists a constant $C>0$, independent of $0<a<a^*$, such that
\begin{equation}\label{zf1}
\|w_a(x)\|_{L^{\infty}(\R^N)}\leq C.
\end{equation}

Actually, for any $L>0$ and $\alpha>0$, set 
\begin{equation}\label{fr1}
\gamma(t)=tt_L^{2\alpha}\,\ \mbox{for}\,\ t>0, \,\ \mbox{where}\,\ t_L=\min\{t,L\}.
\end{equation}
Since $\gamma(t)$ is a nondecreasing function, for any $a,b\geq0$, we have
\begin{equation*}
(b-a)(\gamma(b)-\gamma(a))\geq0.
\end{equation*}
Define
\begin{equation*}
\Gamma(t)=\int_{0}^{t}(\gamma'(\tau))^{\frac{1}{2}}d\tau.
\end{equation*}
For any $a,b\geq0$ with $a<b$, applying Schwartz inequality, we deduce that
\begin{equation}\label{2.7}
\begin{split}
(b-a)(\gamma(b)-\gamma(a))&=(b-a)\int_{a}^{b}\gamma'(t)dt\\
&=(b-a)\int_{a}^{b}(\Gamma'(t))^2dt\\
&\geq \Big(\int_{a}^{b}\Gamma'(t)dt\Big)^2=\big(\Gamma(b)-\Gamma(a)\big)^2.
\end{split}
\end{equation}
Similarly, one can verify that \eqref{2.7} also holds for any $a\geq b$. This further implies that
\begin{equation*}
(b-a)(\gamma(b)-\gamma(a))\geq \big|\Gamma(b)-\Gamma(a)\big|^2\ \ \mbox{for any}\ \ a, b\geq 0.
\end{equation*}
Together with \eqref{fr1}, it then yields that
\begin{equation}\label{02.9}
\big|\Gamma(w_a(x))-\Gamma(w_a(y))\big|^2\leq \big(w_a(x)-w_a(y)\big)\big((w_a(w_a)_{L}^{2\alpha})(x)-(w_a(w_a)_{L}^{2\alpha})(y)\big).
\end{equation}
Taking $\gamma(w_a)=w_a(w_a)_{L}^{2\alpha}\in H^s(\R^N)$ as a test function in \eqref{R1}, we then derive from \eqref{02.9} that
\begin{equation}\label{12.9}
\begin{split}
\big[\Gamma(w_a)\big]_{s}^{2}&\leq \iint_{\R^{2N}}\frac{(w_a(x)-w_a(y))}{|x-y|^{N+2s}}\big((w_a(w_a)_{L}^{2\alpha})(x)-(w_a(w_a)_{L}^{2\alpha})(y)\big)dxdy\\
&\leq\int_{\R^N}c(x)w_a^2(w_a)_{L}^{2\alpha}(x)dx.
\end{split}
\end{equation}

Note that $\Gamma(t)\geq \frac{tt_{L}^{\alpha}}{1+\alpha}$ for any $t\geq0$. Applying Sobolev's inequality, we then deduce that
\begin{equation}\label{c1}
\big[\Gamma(w_a)\big]_s^2\geq S_*\|\Gamma(w_a)\|^{2}_{L^{2^{*}_{s}}(\R^N)}\geq\Big(\frac{1}{1+\alpha}\Big)^{2}S_*\|w_a(w_a)_{L}^{\alpha}\|^{2}_{L^{2^{*}_{s}}(\R^N)}.
\end{equation}
We derive from \eqref{12.9} and \eqref{c1} that
\begin{equation}\label{jc2}
\begin{split}
\Big(\frac{1}{1+\alpha}\Big)^{2}S_*\|w_a(w_a)_{L}^{\alpha}\|^{2}_{L^{2^{*}_{s}}(\R^N)}&\leq \int_{\R^N}c(x)w_a^2(w_a)_{L}^{2\alpha}(x)dx\\
&=\int_{\R^N}c(x)(w_a(w_a)_{L}^{\alpha}(x))^2dx.
\end{split}
\end{equation}
Recall from Lemma \ref{A.1} that $c(x)\!\in\! L^{p}(\R^N)$, where $p\!\in\! \big(\frac{N}{2s},\frac{N^2}{2s(N+b-2s)}\big)$. By H\"{o}lder's inequality, we obtain that
\begin{equation}\label{c3}
\int_{\R^N}c(x)(w_a(w_a)_{L}^{\alpha}(x))^2dx\leq \|c(x)\|_{L^p(\R^N)}\big\|w_a(w_a)_{L}^{\alpha}(x)\big\|^2_{L^{\frac{2p}{p-1}}(\R^N)}.
\end{equation}
Since $p\!\in\! \big(\frac{N}{2s},\frac{N^2}{2s(N+b-2s)}\big)$, it yields that $\frac{2p}{p-1}\in(2,2_s^*)$. By Young's inequality and interpolation inequality, we then deduce from \eqref{jc2} and \eqref{c3} that there exists $\theta\in(0,1)$ such that
\begin{equation}\label{c4}
\begin{split}
&\|w_a(w_a)_{L}^{\alpha}\|^{2}_{L^{2^{*}_{s}}(\R^N)}\\
\leq&C(1+\alpha)^2\varepsilon\|c(x)\|_{L^{p}(\R^N)}\|w_a(w_a)_{L}^{\alpha}(x)\|^2_{L^{2^{*}_{s}}(\R^N)}\\
&+\frac{C}{\varepsilon^{\frac{\theta}{1-\theta}}}(1+\alpha)^2\|c(x)\|_{L^{p}(\R^N)}\|w_a(w_a)_{L}^{\alpha}(x)\|^2_{L^{2}(\R^N)}.
\end{split}
\end{equation}
Taking $\varepsilon>0$ such that
\begin{equation*}
C(1+\alpha)^2\varepsilon\|c(x)\|_{L^{p}(\R^N)}=\frac{1}{2},
\end{equation*}
it then follows from \eqref{c4} that
\begin{equation}\label{jc5}
\begin{split}
\|w_a(w_a)_{L}^{\alpha}\|^{2}_{L^{2^{*}_{s}}(\R^N)}&\leq C(1+\alpha)^{\frac{2}{1-\theta}}\|w_a(w_a)_{L}^{\alpha}(x)\|^2_{L^{2}(\R^N)}\\
&:=M_{\alpha}\|w_a(w_a)_{L}^{\alpha}(x)\|^2_{L^{2}(\R^N)}.
\end{split}
\end{equation}
Letting $L\to\infty$, by the monotone convergence theorem, we then obtain from \eqref{jc5} that
\begin{equation}\label{0c6}
\|w_{a}^{\alpha+1}\|^{2}_{L^{2^{*}_{s}}(\R^N)}\leq M_{\alpha}\|w_{a}^{\alpha+1}\|^{2}_{L^{2}(\R^N)}.
\end{equation}
Note that
\begin{equation*}
M_{\alpha}\leq C(1+\alpha)^{\frac{2}{1-\theta}}\leq M_{0}^{2}e^{2\sqrt{\alpha+1}},
\end{equation*}
where $M_0$ is independent of  $\alpha$. Combining this with \eqref{0c6}, it further indicates that
\begin{equation*}\label{jc7}
\|w_a\|_{L^{2^{*}_{s}(\alpha+1)}(\R^N)}\leq M_{0}^{\frac{1}{1+\alpha}}e^{\frac{1}{\sqrt{1+\alpha}}}\|w_a\|_{L^{2(\alpha+1)}(\R^N)}.
\end{equation*}

Choosing $\alpha_0=0$ and $2(\alpha_{n+1}+1)=2_s^*(\alpha_{n}+1)$, by iteration, we deduce that
\begin{equation}\label{jc8}
\|w_a\|_{L^{2^{*}_{s}(\alpha_n+1)}(\R^N)}\leq M_{0}^{\sum_{i=0}^{n}\frac{1}{1+\alpha_i}}e^{\sum_{i=0}^{n}\frac{1}{\sqrt{1+\alpha_i}}}\|w_a\|_{L^{2(\alpha_0+1)}(\R^N)}.
\end{equation}
Since $1+\alpha_n=\Big(\frac{N}{N-2s}\Big)^n$, we have
\begin{equation*}\label{jc9}
\sum_{i=0}^{\infty}\frac{1}{1+\alpha_i}<\infty\ \ \ \mbox{and}\ \ \ \sum_{i=0}^{\infty}\frac{1}{\sqrt{1+\alpha_i}}<\infty.
\end{equation*}
Together with \eqref{jc8}, it then yields that
\begin{equation}\label{c10}
\|w_a\|_{L^\infty(\R^N)}=\lim\limits_{n\rightarrow\infty}\|w_a\|_{L^{2^{*}_{s}(\alpha_n+1)}(\R^N)}<C\|w_a\|_{L^{2}(\R^N)}=C,
\end{equation}
and \eqref{zf1} is thus proved.

\vskip 0.05truein

\noindent{\em  Step 2.} We prove that
\begin{equation}\label{0c12}
w_a(x)\rightarrow 0 \ \ \mbox{as}\ \ |x|\rightarrow \infty\ \ \mbox{uniformly\ \ for}\ \ a\nearrow a^*.
\end{equation}
Indeed, denote
\begin{equation*}\label{c13}
g_a(x):=a w_a^{1+2\beta^2} |x|^{-b}\ \ \mbox{and}\ \ h_a(x):= g_a(x)-\varepsilon_a^{2s}V(\varepsilon_a x),
\end{equation*}
we deduce from \eqref{2.9} and \eqref{2.10} that
\begin{equation}\label{0c14}
(-\Delta)^{s}w_a(x)+\eta w_a(x)=h_a(x)\ \ \mbox{in}\ \ \R^N,\ \ \mbox{where}\ \ \eta=-\varepsilon_a^{2s}\mu_a >\frac{\beta^2}{2}>0\ \ \mbox{as}\ \ a\nearrow a^*.
\end{equation}
For any $\lambda>0$, let $G_\lambda$ be the Fourier transform of $(|\xi|^{2s}+\lambda)^{-1}$.  We note from \cite [Lemma C.1] {FLS} that there exists a constant $K_1>0$ such that
\begin{equation}\label{c16}
0<G_{\eta}(x)\leq \frac{K_1}{|x|^{N+2s}}\ \ \mbox{for\ \ any}\ \ |x|>1,
\end{equation}
and 
\begin{equation}\label{jc16}
G_{\eta}(x)\in L^q(\R^N)\ \ \mbox{for\ \ any}\ \ q\in \Big[1,\frac{N}{N-2s}\Big).
\end{equation}
 It then follows from \eqref{0c14} that
\begin{equation}\label{0c15}
\begin{split}
0\leq w_a(x)=G_{\eta}\ast h_a(x)&=\int_{\R^N}G_{\eta}(x-y)h_a(y)dy\\
&= \int_{\R^N}G_{\eta}(x-y)g_a(y)dy-\varepsilon_a^{2s} \int_{\R^N}G_{\eta}(x-y)V(\varepsilon_a y)dy\\
&:=B_1-B_2\leq B_1.
\end{split}
\end{equation}

We now estimate $B_1$ as follows. Motivated by \cite{AM, AV}, for any $0<\delta<1$, we define
\begin{equation*}\label{c17}
A_{\delta}:=\big\{y\in\R^N: |y-x|\geq \frac{1}{\delta}\big\}\ \ \ \mbox{and}\ \ \ B_{\delta}:=\big\{y\in\R^N: |y-x|< \frac{1}{\delta}\big\}.
\end{equation*}
We deduce that
\begin{equation}\label{c18}
\begin{split}
0\leq B_1&=\int_{A_{\delta}}G_{\eta}(x-y)g_a(y)dy+\int_{B_{\delta}}G_{\eta}(x-y)g_a(y)dy\\
&:=D_1+D_2.
\end{split}
\end{equation}
Since $w_a\in L^{\infty}(\R^N)$ and $w_a\rightarrow w_0$ in $H^s(\R^N)$ as $a\nearrow a^*$, we derive that  $w_a\rightarrow w_0$ in $L^{r}(\R^N)$ for any $r\in[2,\infty)$. Using a similar argument to \eqref{2.1} then yields that
\begin{equation}\label{c19}
g_a(x)\rightarrow a^*w_0^{1+2\beta^2}|x|^{-b}:=g(x)\ \ \mbox{in}\ \ L^{m}(\R^N)\ \ \mbox{as}\ \ a\nearrow a^*\ \ \mbox{for any}\ \ m\in\Big[1, \frac{N}{b}\Big).
\end{equation}
We thus derive from \eqref{c16} that
\begin{equation}\label{c20}
\begin{split}
D_1&\leq K_1\int_{A_{\delta}}\frac{g_a(y)}{|x-y|^{N+2s}}dy\leq K_1 \|g_a\|_{L^2(A_{\delta})}\Big(\int_{A_{\delta}}\frac{1}{|x-y|^{2N+4s}}dy\Big)^{\frac{1}{2}}\\
&\leq K_1\delta^{2s} \|g_a\|_{L^2(A_{\delta})}\Big(\int_{A_{\delta}}\frac{1}{|x-y|^{2N}}dy\Big)^{\frac{1}{2}}\leq C \delta^{2s}.
\end{split}
\end{equation}

As for $D_2$, we deduce from \eqref{c18} that
\begin{equation}\label{c21}
D_2\leq\int_{B_{\delta}}G_{\eta}(x-y)|g_a(y)-g(y)|dy+\int_{B_{\delta}}G_{\eta}(x-y)g(y)dy.
\end{equation}
Since $\frac{b}{N}+\frac{N-2s}{N}<1$, there exist $m\in [1,\frac{N}{b})$ and $q\in[1,\frac{N}{N-2s})$ such that $\frac{1}{m}+\frac{1}{q}=1$. Applying H\"{o}lder's inequality, one can obtain from \eqref{c21} that
\begin{equation*}\label{c22}
D_2\leq \|G_{\eta}\|_{L^q(\R^N)} \|g_a-g\|_{L^m(\R^N)}+ \|G_{\eta}\|_{L^q(\R^N)} \|g\|_{L^m(B_{\delta})}.
\end{equation*}
Noting that $\|g_a-g\|_{L^m(\R^N)}\rightarrow0$ as $a\nearrow a^*$ and $\|g\|_{L^m(B_{\delta})}\rightarrow0$ as $|x|\rightarrow\infty$, we then deduce from \eqref{jc16} that
\begin{equation*}\label{c23}
D_2\leq \delta\ \ \mbox{as}\ \ |x|\to\infty \ \ \mbox{uniformly\ \ for}\ \  a\nearrow a^*.
\end{equation*}
Together with \eqref{c20}, it further yields that
\begin{equation*}\label{c24}
B_1=\int_{\R^N}G_{\eta}(x-y)g_a(y)dy\leq C \delta^{2s}+\delta\ \ \mbox{as}\ \ |x|\to\infty \ \ \mbox{uniformly for}\ \  a\nearrow a^*.
\end{equation*}
Let $\delta\rightarrow0$, we deduce that $B_1\rightarrow0$ as $|x|\rightarrow\infty$ uniformly for $a\nearrow a^*$. We then derive from \eqref{0c15} that \eqref{0c12} holds true. The proof of Theorem \ref{th3} is therefore complete.\qed

\vskip 0.05truein	
\noindent {\bf Acknowledgements:}  The authors thank Professor Yujin Guo very much for his fruitful discussions on the present paper.

\end{document}